\documentclass[reqno]{article}
\usepackage{tikz}

\usepackage{amsmath,amsthm,amssymb}
\usepackage{graphicx}
\graphicspath{ {./Steiner/} }
\newtheorem{thm}{Theorem}

\newtheorem{conj}[thm]{Conjecture}

\def\beq{\begin{equation}}\def\eeq{\end{equation}}
\def\beqn{\begin{eqnarray}}\def\eeqn{\end{eqnarray}}

\def\qed{\ifhmode\unskip\nobreak\fi\quad\ifmmode\Box\else$\Box$\fi}

\newcommand{\comment}[1]{}

\begin{document}
\title{Large Monochromatic Components of Small Diameter\thanks{2010 Mathematics Subject Classification 05C55. Keywords: Ramsey theory, monochromatic component, diameter. This research was done under the auspices of the Budapest Semesters in Mathematics program}}

\author{Erik Carlson\thanks{Department of Mathematics and Statistics, Carleton College, Northfield, Minnesota, USA, \texttt{egcarlson57@gmail.com}} \and
Ryan R. Martin\thanks{Department of Mathematics, Iowa State University, Ames, Iowa, USA, \texttt{rymartin@iastate.edu}. Research supported in part by a Simons Foundation Collaboration Grant (\#353292) and by the J. William Fulbright Educational Exchange Program. Corresponding author.} \and
Bo Peng\thanks{Department of Mathematics and Statistics, Carleton College, Northfield, Minnesota, USA, \texttt{pengbo1998@yahoo.com}} \and
Mikl\'os Ruszink\'o\thanks{Alfr\'ed R\'enyi Institute of Mathematics, Budapest, P.O. Box 127, Budapest, Hungary, H-1364, Faculty of Information Technology and Bionics, P\'azm\'any P\'eter Catholic University, and NYU Abu Dhabi, P.O. Box 129188 Saadiyat Island, Abu Dhabi, United Arab Emirates \texttt{ruszinko.miklos@renyi.hu}. Research supported in part by NKFIH Grant No. 132696.}}

\maketitle

\begin{abstract}
	Gy\'arf\'as~\cite{Gy2} conjectured in 2011 that every $r$-edge-colored $K_n$ contains a monochromatic component of bounded (``perhaps three'') diameter on at least $n/(r-1)$ vertices. 
	Letzter~\cite{L} proved this conjecture with diameter four. 
	In this note we improve the result in the case of $r=3$: We show that in every $3$-edge-coloring of $K_n$ either there is a monochromatic component of diameter at most three on at least $n/2$ vertices or \textit{every} color class is spanning and has diameter at most four.
\end{abstract}

\section{Monochromatic components}
	An easy exercise in an introductory graph theory course -- a remark by Erd\H{o}s and Rado, see~\cite{Gy2} -- states that any $2$-coloring of the edges of $K_n$ has a monochromatic spanning component. 
	In general, Gy\'{a}rf\'{a}s \cite{Gy1} proved that the largest monochromatic component in an $r$-edge-coloring of $K_n$ has order at least $n/(r-1)$ and equality holds if an affine plane of order $r-1$ exists and $(r-1)^2$ divides $n$. 
	F\"uredi \cite{F} proved the significantly larger lower bound $n/(r-1-(r-1)^{-1})$ in the case that there is no affine plane of order $r-1$. 
	This connection to the existence of affine planes suggests that to determine exactly the maximum size of a monochromatic component is extremely difficult in general.

	A double [triple] star is the tree obtained by joining the centers of two [three] stars by  a path of length one [two].  
	Clearly, double or triple stars have diameter 3 or 4 respectively. 
	Additional structure on large monochromatic components has been conjectured by Gy\'arf\'as~\cite{Gy2}.

	\begin{conj}[Gy\'arf\'as, Problem 4.2 in \cite{Gy2}] \label{P1}
		For $r\ge 3$, is there a monochromatic double star on at least $n/(r-1)-o(n)$ vertices in every $r$-coloring of $K_n$?
	\end{conj}

	A weaker version of the problem reads as follows.

	\begin{conj}[Gy\'arf\'as, Problem 4.3 in \cite{Gy2}] \label{P}
		Given positive numbers $n$, $r$.
		Is there a constant $d$ (perhaps $d=3$) such that in every $r$-coloring of $K_n$ there is a monochromatic subgraph of diameter at most $d$ with at least $n/(r-1)$ vertices?
	\end{conj}

	The assumption $r\ge 3$ in Conjecture \ref{P1} is necessary, since a random two-coloring will a give a monochromatic double star of size $\approx 3n/4$ only. The best result for double stars is due to Gy\'arf\'as and S\'ark\"ozy.

	\begin{thm}[Gy\'arf\'as, S\'ark\"ozy~\cite{GS}] \label{GS}
		Every $r$-edge-coloring of $K_n$ contains a monochromatic double star on at least $\frac{n(r+1)+r-1}{r^2}$ vertices.
	\end{thm}

	The bipartite Ramsey number of the double star has been determined by Mubayi \cite{M}.
	The result of Theorem~\ref{bip} is tight if each color class is biregular. 
	\begin{thm}[Mubayi~\cite{M}]\label{bip} 
		In every $r$-edge-coloring of the complete bipartite graph $K_{k,\ell }$ there is a monochromatic double star of order $\frac{k+\ell}{r}$.
	\end{thm}

	The weaker Conjecture~\ref{P} was later shown to be true by Ruszink\'o~\cite{R} with $d=5$.
	\begin{thm}[Ruszink\'o~\cite{R}] \label{fo1}
		In every $r$-edge-coloring of $K_n$ there is a monochromatic subgraph of diameter at most $5$ on at least $n/(r-1)$ vertices.
	\end{thm}

	This was further improved and was shown to be true for $d=4$ by Letzter~\cite{L}.
	\begin{thm}[Letzter~\cite{L}] \label{fo2}
		In every $r$-edge-coloring of $K_n$ there is a monochromatic triple star on at least $n/(r-1)$ vertices.
	\end{thm}

	For the case of $d=r=2$, the following tight bound was proved by Erd\H{o}s and Fowler~\cite{EF}.
	\begin{thm}[Erd\H{o}s, Fowler \cite{EF}] \label{EF1}
		Every $2$-edge-coloring of $K_n$ contains a monochromatic connected subgraph of diameter at most $2$ on at least $3n/4$ vertices.
	\end{thm}
	Moreover, for $r=3,4,5,6$, Ruszink\'o, Song and Szabo~\cite{RSS} constructed colorings where the maximum size of a monochromatic, diameter 2 subgraph is strictly less than $n/(r-1)$, suggesting that $d=3$ is best possible for diameter in Conjecture~\ref{P}.\\

	In this note we further improve (in terms of diameter) Theorem \ref{fo2} for three colors. 
	Let $G_\alpha$, $G_\beta$, and $G_\gamma$ be the subgraphs of $K_n$ induced by the edges that have color $\alpha$, $\beta$, and $\gamma$ respectively.

	\begin{thm}\label{main}
		In every $3$-edge-coloring of $K_n$ either there is a monochromatic connected subgraph of diameter at most $3$ on at least $n/2$ vertices or each of $G_\alpha$, $G_\beta$, and $G_\gamma$ is spanning and has diameter at most $4$.
	\end{thm}

	\begin{proof}
		By Theorem~\ref{bip}, we may assume that each of $G_\alpha$, $G_\beta$, and $G_\gamma$ is both spanning and connected because if one is not, then the union of the other two color classes is a complete bipartite graph on $n$ vertices.
		
		Suppose, towards a contradiction and without loss of generality, that the distance between $w_1$ and $w_2$ is at least $5$ in $G_\alpha$ and $w_1w_2\in E(G_\beta)$.
		The set of the vertices $U$ of the double star centered by $w_1$ and $w_2$ in $G_\beta$ must contain less than $n/2$ vertices, otherwise the theorem is proven.
		
		Note that there are no $\beta$-colored edges from $\{w_1,w_2\}$ to $V\setminus U$ by definition.
		Split the remaining vertices of $V\setminus U$ into $3$ parts: 
		\begin{align*}
			X 	&=	\{v\in V\setminus U:~vw_1\in E(G_\gamma),~vw_2\in E(G_\alpha)\} , \\
			Y 	&= 	\{v\in V\setminus U:~vw_1\in E(G_\gamma),~vw_2\in E(G_\gamma)\} , \\
			Z 	&= 	\{v\in V\setminus U:~vw_1\in E(G_\alpha),~vw_2\in E(G_\gamma)\} . 
		\end{align*}
		Note that there are no vertices $v$, such that $vw_1\in E(G_\alpha)$ and $vw_2\in E(G_\alpha)$, or else the distance between $w_1$ and $w_2$ in $G_1$ would be $2<5$. 
		Clearly, neither $X$ nor $Z$ is empty, or else there is a star in $G_\gamma$ (centered at either $w_1$ or $w_2$) of order greater than $n/2$. 
		Furthermore, no edge between $X$ and $Z$ is colored $\alpha$, or else we have a path of length $3$ in $G_\alpha$ between $w_1$ and $w_2$.
		
		In addition, there is a length $2$ path in color $\gamma$ between each pair of vertices in $X$ (through $w_1$), between each vertex in $X$ with each vertex in $Y$ (through $w_1$), between each pair of vertices in $Y$ (through either $w_1$ or $w_2$), between each vertex in $Y$ and each vertex in $Z$ (through $w_2$), and between each pair of vertices in $Z$ (through $w_2$). 
		
		Since $X\cup Y\cup Z\cup\{w_1,w_2\}$ contains more than $n/2$ vertices, there must exist some vertices, $v_X\in X$, $v_Z\in Z$, such that their distance in color $\gamma$ within the vertex set $X\cup Y\cup Z\cup\{w_1,w_2\}$ is at least $4$, otherwise we have found a vertex set of diameter $3$ in color $\gamma$ of size larger than $n/2$.
		For this to be the case, neither $v_X$ nor $v_Z$ can have an edge colored $\gamma$ connecting it to any vertex in $Y$, otherwise there would be a path of length 3 connecting $v_X$ and $v_Z$ in $G_\gamma$. 
		Furthermore, since there is no edge of color $\alpha$ between $X$ and $Z$, $v_X$ must have only edges colored $\beta$ between itself and all vertices of $Z$, and $v_Z$ must have only edges colored $\beta$ between itself and all vertices of $X$, otherwise, again $v_X$ and $v_Z$ would have distance at most 3 in $G_\gamma$. 

		Now, we have a double star in color $\beta$, anchored at $v_X$ and $v_Z$, containing all of $X\cup Z$. 
		If $Y$ is empty, the theorem is proven.
		Therefore, there must be some $v_Y \in Y$ such that neither $v_Yv_X$ nor $v_Yv_Z$ has color $\beta$, otherwise there is a double star in color $\beta$, anchored at $v_X$ and $v_Z$, containing $X\cup Y\cup Z$, which is a double star on at least $n/2$ vertices.

		So, the edges $v_Yv_X$ and $v_Yv_Z$ have neither color $\beta$ nor color $\gamma$. 
		Therefore, both such edges have color $\alpha$. 
		This produces a path in $G_\alpha$, namely $w_1v_Zv_Yv_Xw_2$, of length $4$, which is a contradiction to the assumption that $w_1$ and $w_2$ have distance at least $5$ in $G_\alpha$. 	
	\end{proof}
	
\medskip
\noindent
{\bf Conclusion}. Though Theorem~\ref{main} does not prove Conjecture~\ref{P} for $d=3$ in the case of three colors, it gives support to this very natural and surprisingly difficult question. 

\medskip
\noindent
{\bf Acknowledgements}. This research was part of a class in the Budapest Semesters in Mathematics program in the Fall of 2019. The authors also wish to acknowledge the R\'enyi Institute of Mathematics for the use of its facilities. The authors would like also to thank G\'abor S\'ark\"ozy and Andr\'as Gy\'arf\'as for fruitful comments and discussions.

\end{document}